\documentclass[letterpaper,oneside,11pt]{article}

\usepackage[USenglish]{babel} 
\usepackage[T1]{fontenc}
\usepackage[ansinew]{inputenc}

\usepackage{lmodern} 

\usepackage{graphicx} 

\usepackage{amsmath}
\usepackage{amsthm}
\usepackage{amsfonts}

\begin{document}

\pagestyle{empty} 


\title{On Eisenstein's formula for the Fermat quotient}
\author{John Blythe Dobson (j.dobson@uwinnipeg.ca)}

\maketitle

\tableofcontents 
\cleardoublepage 

\pagestyle{plain} 


\begin{abstract}
\noindent
This paper presents some refinements of the representation of the Fermat quotient of base 2 as an alternating series which was discovered by Eisenstein in 1850, including some evaluations that are believed to be new.

\end{abstract}

\section{Introduction}
For the Fermat quotient $q_p(b) = (b^{p - 1} - 1)/p \pmod{p}$ we employ wherever possible the briefer notation $q_b$. Henceforth all congruences are assumed to be mod $p$ unless otherwise stated, and $\lfloor \cdot \rfloor$ signifies the greatest-integer function.

The fact that Fermat quotients can be expressed as sums involving reciprocals of integers in $\left\{ 1, p-1\right\}$ was discovered in 1850 for the case $b = 2$ by Eisenstein (\cite{Eisenstein}, p.\ 41), who gives

\begin{equation} \label{eq:Eisenstein}
2 \cdot q_2 \equiv 1 - \frac{1}{2} + \frac{1}{3} - \frac{1}{4} + \ldots + \frac{1}{(p - 2)} - \frac{1}{(p - 1)} = \sum_{k=1}^{p-1} \frac{(-1)^{k-1}}{k}.
\end{equation}

\noindent
Subsequent researches in this direction obtained comparable results entailing fewer terms. Stern (\cite{Stern}, p.\ 185) found

\begin{equation} \label{eq:EisensteinHalf}
q_2 \equiv \sum_{k=1}^{(p-1)/2} \frac{(-1)^{k-1}}{k},
\end{equation}

\noindent
and drawing on work of Emma Lehmer \cite{Lehmer}, Zhi-Hong Sun \cite{ZHSun1992}, pt. 1, Corollary 1.3, found

\begin{equation} \label{eq:EisensteinThird}
2 \cdot q_2 \equiv \sum_{k=1}^{\lfloor p/3 \rfloor} \frac{(-1)^{k-1}}{k}.
\end{equation}

\noindent
Subtracting (\ref{eq:EisensteinThird}) from (\ref{eq:EisensteinHalf}) gives

\begin{equation} \label{eq:EisensteinSixth}
-q_2 \equiv \sum_{k=\lfloor p/3 \rfloor}^{(p-1)/2} \frac{(-1)^{k-1}}{k}.
\end{equation}

\noindent
Dilcher \& Skula (\cite{Dilcher}, p.\ 389) report numerical investigations of all possible instances of $s(k, N)$ with $N \le 46$ for the two values of $p$ for which $q_2$ vanishes, namely 1093 and 3511. We have extended these calculations to consider all possible cases of $N$, without finding any new cases of such sums that consistently vanish with $q_2$. This constitutes a numerical proof that (\ref{eq:EisensteinSixth}) is optimal. The demonstration by Wieferich and Mirimanoff of the connection between Fermat quotients and the first case of Fermat's Last Theorem (FLT) gave such results fresh interest, which they retain despite the full proof of FLT by Wiles in 1995.

While sums like Eisenstein's (\ref{eq:Eisenstein}) lose their connection with the Fermat quotient when their range includes less that $\frac{1}{6}$ of the terms in the interval $\left\{ 1, p-1\right\}$, nevertheless they possess arithmetical interest of their own; and a new motivation for their study was supplied by the landmark 1992 paper of Dilcher \& Skula just mentioned, which showed that a failure of the first case of FLT would require the simultaneous vanishing of closely related sums on a massive scale, a point to which we shall return below.

In an earlier paper \cite{Dobson}, we have reviewed previous researches pertaining to the Fermat quotient, and presented a few new ones, deliberately avoiding the use of alternating series in order to faciliate comparison with the results of Stern \cite{Stern}, Lerch \cite{Lerch}, Skula (\cite{Skula1992}, \cite{Skula2008}), and those who have followed their conventions. In the present paper, however, we return to the original inspiration of Eisenstein's formula, and explore results of independent interest which find their most natural expression as alternating series. The numerous works of Zhi-Hong Sun in this area furnish some particularly enigmatic examples of such series, which are not adequately explained by the presentation in \cite{Dobson}, and it seemed that a fresh approach was needed. The most conspicuous of these examples is Sun's determination in \cite{ZHSun1992}, pt.\ 3, Theorem 3.2, nos.\ 3 and 5, that

\begin{equation} \label{eq:Sun15}
\sum_{k=\lfloor p/5 \rfloor}^{\lfloor p/3 \rfloor} \frac{(-1)^k}{k} \equiv -\sum_{k=1}^{\lfloor 2p/5 \rfloor} \frac{(-1)^k}{k}.
\end{equation}

\noindent
That the summand should have a corresponding (if opposite) evaluation for these choices of range is surely surprising --- note that the sum on the right runs over only $\frac{2}{15}$ of the range $\left\{1, p-1\right\}$ --- but it will be demonstrated that such results can be proved in a natural and uniform manner, independent of a consideration of Fibonacci quotients and other entities that supplied the original motivation for their discovery.

Here, as in the previous paper, we use the notations

\begin{displaymath}
s(k, N) = \sum_{\substack{j=\lfloor\frac{kp}{N}\rfloor + 1\\ j \neq p}}^{\lfloor\frac{(k + 1)p}{N}\rfloor}
\frac{1}{j} \textrm{\hspace{1em}\hspace{0.5em} and \hspace{1em}\hspace{0.5em}} s^{\ast}(k, N) = \sum_{\substack{j=\lfloor\frac{kp}{N}\rfloor + 1\\ j \neq p}}^{\lfloor\frac{(k + 1)p}{N}\rfloor} \frac{(-1)^j}{j},
\end{displaymath}

\noindent
and it is always assumed that $p$ is sufficiently large that $s(k, N)$ contains at least one element; the provision $j \neq p$ is necessary in when $k+1=N$. At the risk of belaboring the obvious, we point out the equivalency of the following relations, which can be useful in recognizing values of $s^{\ast}(k, N)$ buried in the literature:

\begin{equation} \label{eq:ThreeRules}
\begin{split}
s^{\ast}(k, N) & = \hphantom{-}s{''}(k, N) - s^\prime(k, N) \\
               & = \hphantom{-}s(k, N) - 2s^\prime(k, N) \\
               & = -s(k, N) + 2s{''}(k, N) \\
\end{split}
\end{equation}

\noindent
where $s^\prime(k, N)$ and $s{''}(k, N)$ represent the odd and even terms in $s(k, N)$, respectively. For even $N$, the relation

\begin{equation} \label{eq:FourthRule}
s{''}(k, N) = \frac{1}{2}s(k, \frac{N}{2})
\end{equation}

\noindent
can likewise be used if $s(k, \frac{N}{2})$ happens to be known.

We recall the fact that for $s(k, N)$, the terms in $\left\{\frac{(p+1)}{2}, p-1\right\}$ (descending) are the complements $\pmod{p}$ of those in $\left\{ 1, \frac{(p-1)}{2}\right\}$ (ascending), implying that

\begin{equation}
s(k, N) \equiv -s(N-1-k, N),
\end{equation}

\noindent
while for $s^{\ast}(k, N)$, the corresponding terms are congruent rather than complementary, so that

\begin{equation}
s^{\ast}(k, N) \equiv s^{\ast}(N-1-k, N).
\end{equation}

\noindent
For ease of comparison with previous literature, in our final results we usually restrict $k$ so as to be less than $\frac{p-1}{2}$, but in the proofs we use whichever form seems more intelligible or expressive in the given situation. In the next few paragraphs we shall make frequent use, without further comment, of the following corollary of (\ref{eq:sStar}) below:

\begin{displaymath}
s^{\ast}(0, N) \equiv -s^{\ast}(1, 2N).
\end{displaymath}

Some specific values of $s^{\ast}(k, N)$ have long been known: $s^{\ast}(0, 1) \equiv -s(1, 2) \equiv -2 \cdot q_2$ is the classic result from 1850 of Eisenstein already mentioned, and $s^{\ast}(0, 2) \equiv -s(1, 4) \equiv -q_2$ is derived in an important but often overlooked paper of 1887 by Stern (\cite{Stern}, pp.\ 185--86). Glaisher's statement of the latter result (\cite{Glaisher}, pp.\ 22--23) is a rediscovery and does not warrant the level of credit assigned to it in Dickson's \textit{History} \cite{Dickson}, 1:111. There followed a nearly century-long lull in progress on this problem until 1982, when H.\,C.\ Williams (\cite{Williams1982}, p.\ 369) in effect evaluated $s^{\ast}(0, 5)$, equivalent to $-s(1,10)$ and apparently the earliest result for any value of $s(k, 10)$, in terms of the Fibonacci quotient; it is trivial to deduce from his Theorem 2 that

\begin{displaymath}
s^{\ast}(0, 5) \equiv - 2 \cdot q_2 - \frac{5}{2}\frac{F_{p-\epsilon}}{p},
\end{displaymath}

\noindent
where $\epsilon = \left(\frac{5}{p}\right)$, $\left(\frac{\hphantom{5}}{\hphantom{p}}\right)$ being the Legendre symbol. Then, Zhi-Hong Sun, in his monumental contribution of 1992--1993 \cite{ZHSun1992}, evaluated $s^{\ast}(0, 3) \equiv -2 \cdot q_2$ (Corollary 1.3); $s^{\ast}(1, 3) \equiv 2 \cdot q_2$ (implicitly, in Corollary 1.4); $s^{\ast}(0, 4)$, equivalent to $-s(1, 8)$, in terms of $q_2$ and the Pell sequence (Lemma 2.2); $s^{\ast}(1, 4)$ (Theorem 3.4, no.\ iii); $s^{\ast}(1, 5)$ (implicitly, in Theorem 3.2, no.\ 3); and $s^{\ast}(3, 12)$, in terms of $q_2$ and the Pell sequence, apparently without recognizing that it is equivalent to $s(2, 8)$ (Theorem 3.4, no. iv). As we have seen (\ref{eq:EisensteinSixth}), his argument leads easily to the evaluation $s^{\ast}(2, 6) \equiv q_2.$ Most remarkably, he was able to evaluate $s^{\ast}(0, 9)$, equivalent to $-s(1, 18)$ and still apparently the only known result for any $s(k, 18)$, in terms of $q_2$ and a complex recurrence relation (Corollary 2.4). He also gave a new expression for $s^{\ast}(0, 5)$ in terms of the Lucas quotient (Corollary 1.11, no.\ 1):

\begin{displaymath}
s^{\ast}(0, 5) \equiv -2 \cdot q_2 - \frac{\left(\frac{5}{p}\right)L_{p+\epsilon}-3}{p},
\end{displaymath}

\noindent
where as before $\epsilon = \left(\frac{5}{p}\right)$. Improvements on some of the classical results may be obtained from another paper of the same author \cite{ZHSun2008} by use of the rules given in (\ref{eq:ThreeRules}) and (\ref{eq:FourthRule}). Sun there gives an evaluation of $s^\prime(0, 1)$ mod $p^3$, and subtracting twice this quantity from the well-known evaluation $s(0, 1) \equiv \frac{1}{3}p^2 \cdot B_{p-3} \pmod{p^3}$, where $B$ is a Bernoulli number, gives

\begin{displaymath}
s^{\ast}(0, 1) \equiv-2 \cdot q_2 + p \cdot q_2^2 - \frac{2}{3}p^2 \cdot q_2^3 - \frac{1}{4}p^2 \cdot B_{p-3} \pmod{p^3}.
\end{displaymath}

\noindent
Further, since $s^{\ast}(0, 2) \equiv s(0, 4) - s(0, 2)$ and Sun evaluates each of the latter values mod $p^3$, we could write down a similar result for $s^{\ast}(0, 2)$, but it would be rather complicated. Zhi-Wei Sun (\cite{ZWSun2002}, p.\ 17), evaluates $s^{\ast}(0, 6)$, equivalent to $-s(1, 12)$, in terms of the Lucas sequence.

In this brief overview, we have deliberately passed over examples where $s^{\ast}(k, N)$ figures in evaluations of more complex or less well-studied series. Values that can be expressed purely in terms of the Fermat quotient are collected for ease of reference in Table \ref{Table_1}.

\section{Eisenstein's congruence and Fermat's Last Theorem}

Historically, the main motivation for the study of $s^{\ast}(k, N)$ has been its connection with the first case of FLT, either via its relationship with the Fermat quotient or in its own right. Combining results of Dilcher \& Skula \cite{Dilcher}, Cik\'{a}nek \cite{Cikanek}, and the present writer \cite{Dobson}, a failure of the first case of FLT would require $s(k, N) \equiv 0$ for all $N \leq 46$, all oddly even $N \leq 90$, and all remaining $N \leq 94$ for sufficiently large $p$; i.\,e., $p> 5^{(N-1)^{2} (N-2)^{2}/4}$. Thus, to show that even one of these values cannot vanish for a given $p$ proves the first case of FLT for $p$. As pointed out by Dilcher \& Skula (\cite{Dilcher}, p.\ 390, Proposition 10.4), a sum cannot vanish unless it contains at least three terms, since the reciprocals in the summand are distinct mod $p$. This consideration applies regardless of whether or not the signs of the terms alternate (since subtracting a value from a sum is tantamount to adding its complement), so that the first case of FLT follows for $p < 3N$. If we can find sums with higher values of $N$ (and thus smaller ranges) that can be expressed as linear combinations of sums with smaller $N$, we can extend the conditions on any hypothetical counterexample to the first case of FLT.

\section{The evaluation of $s^{\ast}(k, N)$}

\noindent
At first glance, this might not seem to be an especially promising prospect in the case of the alternating sums $s^{\ast}(k, N)$. We recall some formulae from our earlier paper:

\begin{equation} \label{eq:sStar}
s^{\ast}(k, N) \equiv \frac{1}{2} \cdot s(k, 2N) - \frac{1}{2} \cdot s(N + k, 2N);\\
\end{equation}

\begin{equation} \label{eq:sStar2}
s^{\ast}(k, N) + s(k, N) \equiv s(k, 2N);
\end{equation}

\begin{equation} \label{eq:sStar3}
s^{\ast}(k, N) - s(k, N) \equiv s(N-1-k, 2N).
\end{equation}

\noindent
In the following section, we shall need the last in the more homogeneous form

\begin{equation} \label{eq:starSplitInThree}
s^{\ast}(k, N) \equiv s(N-1-k, 2N) + s(2k, 2N) + s(2k+1, 2N).
\end{equation}

\noindent
Note that all of these formulae relate $s^{\ast}(k, N)$ to combinations of sums at least one of which has a \textit{higher} value of $N$, contrary to our goal. And this difficulty persists even in most cases where the above formulae can be simplified, such as the first of the following two:







\subsection{$k=(N-1)/2$ ($N$ odd)}

\noindent
For odd $N$,

\begin{equation} \label{eq:Odd}
s^{\ast}\left(\frac{N-1}{2}, N\right) \equiv s\left(\frac{N-1}{2}, 2N\right).
\end{equation}

\subsection{$k=(N-2)/2$ ($N$ oddly even)}

\noindent
When $N \equiv 2 \pmod{4}$,

\begin{equation} \label{eq:OddlyEven}
s^{\ast}\left(\frac{N-2}{2}, N\right) \equiv \frac{1}{2} \cdot s\left(\frac{N-2}{4}, N\right).
\end{equation}

\noindent
Only in this latter case does the corresponding evaluation in terms of $s(k, N)$ not entail a higher value of $N$.

\section{Some special values of $s^{\ast}(k, N)$}

\noindent
We shall now show that the basic rules gathered in the previous section do not exhaust all the values of $s^{\ast}(k, N)$ that can be evaluated in terms of $s(k, N)$. The surprising relation (\ref{eq:Sun15}) may be written in our notation as

\begin{equation} \label{eq:Sun15rewritten}
s^{\ast}(0, 3) - s^{\ast}(0, 5) \equiv -s^{\ast}(0, 5) - s^{\ast}(1, 5).
\end{equation}

\noindent
This has two separate and interesting implications. The first is the more obvious inference that

\begin{equation} \label{eq:3against5}
s^{\ast}(0, 3) \equiv - s^{\ast}(1, 5).
\end{equation}

\noindent
In our earlier paper (\cite{Dobson}, Theorem 2) we explained this strange-looking correspondence as a consequence of the sole instance where sums of terms $s(k, N)$ with equally-spaced values of $k$ and different values of $N$, neither dividing the other, each contain precisely two terms:

\begin{equation} \label{eq:Sun15rewritten2}
s(0, 6) + s(2, 6) \equiv -s(1, 10) - s(3, 10) \equiv -4 \cdot q_2.
\end{equation}

\noindent
The second consequence of (\ref{eq:Sun15rewritten}), overlooked in our earlier paper, is obtained by observing that in (\ref{eq:Sun15rewritten}), the left-hand side is by definition equivalent to $s^{\ast}(1, 5) - s^{\ast}(5, 15)$, while the right-hand side may be reduced to $s(1, 5)$ using (\ref{eq:sStar2}). Rearranging,

\begin{equation} \label{eq:Sun15rewritten3}
s^{\ast}(5, 15) \equiv s^{\ast}(1, 5) - s(1, 5) \equiv s(1, 10) - 2s(2, 10) - 2s(3, 10).
\end{equation}

\noindent
To reduce the far right-hand side, we employ two congruences discovered by Skula (\cite{Skula2008}, Theorem 3.2):

\begin{displaymath}
2s(0, 10) + 3s(1, 10) + 2s(2, 10) + 3s(3, 10) + 2s(4, 10) \equiv 0;
\end{displaymath}

\begin{displaymath}
s(0, 10) + 2s(1, 10) + s(4, 10) \equiv 0.
\end{displaymath}

\noindent
Adding the first of these to the far right-hand side of (\ref{eq:Sun15rewritten3}), and subtracting twice the second, gives the much simpler expression

\begin{equation} \label{eq:Sun15rewrittenFinal}
s^{\ast}(5, 15) \equiv s(3, 10),
\end{equation}

\noindent
which finally provides a development of $s^{\ast}(5, 15)$ in terms of $s(k, N)$ with $N$ less than 15.

A major motivation of the present paper was to find a less \textit{ad hoc} method of discovering such appealing congruences as (\ref{eq:Sun15rewrittenFinal}). After eliminating the cases satisfactorily accounted for above, an extensive search was made in the literature for other unexpectedly simple valuations, and the following parallels were found:

\begin{equation} \label{2_6}
s^{\ast}(2, 6) \equiv s(1, 4);
\end{equation}

\begin{equation} \label{3_12}
s^{\ast}(3, 12) \equiv s(2, 8).
\end{equation}

\noindent
The first of these is verifiable, with some effort, from Williams \cite{Williams1991}, p.\ 440, while the second can be worked out from Zhi-Hong Sun \cite{ZHSun1992}, pt.\ 3, Theorem 3.4, no.\ iv. The three relations (\ref{2_6}), (\ref{3_12}), and (\ref{eq:Sun15rewrittenFinal}) can be recognized as cases of (\ref{eq:starSplitInThree}) in which the three terms in the right-hand side are adjacent, with the least term divisible by 3, enabling a reduction of the three terms to a single term of the form $s(\_\_, \frac{2N}{3})$. If $2k$ is the least of the three terms, then $N-1-k = 2k+2$ and $3|k$, so $N = 3k+3 \equiv 3 \pmod{9}$. If on the other hand $N-1-k$ is the least term, then $N-1-k = 2k-1$ and $2k \equiv 1 \Rightarrow k \equiv 2 \pmod{3}$, so $N = 3k \equiv 6 \pmod{9}$. Thus these special values of $s^{\ast}(k, N)$ fall into two infinite families, given by the following congruences (with $x = 0, 1, 2, \dots$):

\begin{subequations} \label{TwoThirds}
\begin{gather}
s^{\ast}(3x, 9x+3)   \equiv s(2x, 6x+2) \\
s^{\ast}(3x+2, 9x+6) \equiv s(2x+1, 6x+4).
\end{gather}
\end{subequations}

\noindent
The first congruence supplies valuations of $s^{\ast}(0, 3)$, $s^{\ast}(3, 12)$, $s^{\ast}(6, 21)$, etc., while the second supplies valuations of $s^{\ast}(2, 6)$, $s^{\ast}(5, 15)$, $s^{\ast}(8, 24)$, etc. They therefore account for all of the ostensibly anomalous relationships found in the literature, besides incidentally supplying an alternate derivation of $s^{\ast}(0, 3)$.





\section{The effect of the vanishing of $s(\_, N)$ or $s^{\ast}(\_, N)$ on other sums}

As mentioned in our previous paper \cite{Dobson}, for odd $N$, the vanishing of $s(k, N)$ for every value of $k$ implies the vanishing of $s(k, 2N)$ for every value of $k$. As is clear from its definition (\ref{eq:sStar}), this would likewise imply the vanishing of $s^{\ast}(k, N)$ for every value of $k$.

It should finally be noted that two of the results recorded in this paper, (\ref{eq:OddlyEven}) and (\ref{TwoThirds}), may in some cases permit the conclusion that the vanishing of both $s^{\ast}(k, N)$ and $s(k, N)$ for given $N$ and $k$ is a necessary condition for the failure of the first case of FLT. In such cases, by virtue of (\ref{eq:sStar2}) and (\ref{eq:sStar3}), the same would hold true for $s(k, 2N)$ and $s(N-1-k, 2N)$.

\begin{table} [hb]
\begin{center}
\caption{Presumably complete list of the sums $s^{\ast}(k, N)$ (with $k < N/2$) that can be evaluated solely in terms of Fermat quotients}
\label{Table_1}
\begin{tabular}{@{}lr@{}}
$s^{\ast}(0, 1)$  & $-2 \cdot q_2$ \\
$s^{\ast}(0, 2)$  & $-q_2$ \\
$s^{\ast}(0, 3)$  & $-2 \cdot q_2$ \\ 
$s^{\ast}(1, 3)$  & $2 \cdot q_2$ \\ 
$s^{\ast}(1, 5)$  & $2 \cdot q_2$ \\ 
$s^{\ast}(2, 6)$  & $q_2$ \\ 
\end{tabular}
\end{center}
\end{table}


\clearpage
\begin{thebibliography}{0}

\bibitem{Cikanek}
Petr Cik\'{a}nek, ``A special extension of Wieferich's criterion,'' Math. Comp. {\bf 62} (1994) 923--930.

\bibitem{Dickson}
Leonard Eugene Dickson, \textit{History of the Theory of Numbers}, 3 vols. New York, 1919.

\bibitem{Dilcher}
Karl Dilcher \& Ladislav Skula, ``A New Criterion for the First Case of Fermat's Last Theorem,'' Math. Comp. {\bf 64} (1995) 363--392.

\bibitem{Dobson}
John Blythe Dobson, ``On Lerch's formula for the Fermat quotient,'' available at http://arxiv.org/abs/1103.3907.

\bibitem{Eisenstein}
[G.] Eisenstein, ``Eine neue Gattung zahlentheoretischer Funktionen, welche von zwei Elementen abh\"{a}ngen und durch gewisse lineare Funktional-Gleichungen definirt werden,'' Berichte K\"{o}nigl. Preu\ss{}. Akad. Wiss. Berlin {\bf 15} (1850) 36--42.

\bibitem{Glaisher}
J.\,W.\,L.\ Glaisher, ``On the Residues of $r^{p-1}$ to Modulus $p^2$, $p^3$, etc.,'' Q. J. Math. Oxford {\bf 32} (1900--1901) 1--27.

\bibitem{Lehmer}
Emma Lehmer. ``On Congruences involving Bernoulli Numbers and the Quotients of Fermat and Wilson,'' Ann. of Math. {\bf 39} (1938) 350--360.

\bibitem{Lerch}
M.\ Lerch, ``Zur Theorie des Fermatschen Quotienten\ldots,'' Math. Ann. {\bf 60} (1905) 471--490.

\bibitem{Skula1992}
Ladislav Skula, ``Fermat's Last Theorem and the Fermat Quotients,'' Comment. Math. Univ. St. Pauli {\bf 41} (1992) 35--54.

\bibitem{Skula2008}
Ladislav Skula, ``A note on some relations among special sums of reciprocals modulo $p$,'' Math. Slovaca {\bf 58} (2008) 5--10.

\bibitem{Stern}
M.\ Stern, ``Einige Bemerkungen \"{u}ber die Congruenz $\frac{(r^{p} - r)}{p} \equiv a \pmod{p}$,'' J. Reine Angew. Math. {\bf 100} (1887) 182--188.

\bibitem{ZHSun1992}
Zhi-Hong Sun, ``[The] Combinatorial Sum $\sum_{k \equiv r \pmod{m}} \binom{n}{k}$ and its Applications in Number Theory'' (in Chinese), J. Nanjing Univ. Math. Biquarterly {\bf 9} (1992): 227--240, {\bf 10} (1993): 205--118, {\bf 12} (1995): 90-102. A very full summary in English is available on the author's website, at http://www.hytc.cn/xsjl/szh/coms1.pdf, coms2.pdf, coms3.pdf.

\bibitem{ZHSun2008}
Zhi-Hong Sun, ``Congruences involving Bernoulli and Euler numbers,'' J. Number Theory {\bf 128} (2008) 280--312.

\bibitem{ZWSun2002}
Zhi-Wei Sun, ``On the Sum $\sum_{k \equiv r \pmod{m}} \binom{n}{k}$ and Related Congruences,'' Israel J. Math. {\bf 128} (2002): 135--56.

\bibitem{Williams1982}
H.\,C.\ Williams, ``A Note on the Fibonacci Quotient\ldots,'' Canad. Math. Bull. {\bf 25} (1982) 366--370.

\bibitem{Williams1991}
H.\,C.\ Williams, ``Some formulas concerning the fundamental unit of a real quadratic field,'' Discrete Math. {\bf 92} (1991) 431--440.

\end{thebibliography}
\end{document}